\documentclass[12pt]{amsart}

 \usepackage{amsmath}
\usepackage{amssymb,latexsym,cite}
\usepackage{amscd}
\usepackage{comment}
 \usepackage{xspace}
 \usepackage{verbatim}
 \usepackage{todonotes}

\newtheorem{theorem}{Theorem}[section]

\newtheorem{corollary}[theorem]{Corollary}
\newtheorem{lemma}[theorem]{Lemma}
\theoremstyle{definition}
\newtheorem{remark}[theorem]{Remark}
\newtheorem{definition}[theorem]{Definition}
\numberwithin{equation}{section}

\newcommand{\tr}{\mathrm{tr}^*}

\newcommand{\chr}[1]{\mathbf{1}\ind{#1}}

\newcommand{\BSA}{\begin{subarray}}
\newcommand{\ESA}{\end{subarray}}
\newcommand{\BAL}{\begin{aligned}}
\newcommand{\EAL}{\end{aligned}}

\newcommand{\note}[1]{\noindent\textit{#1.}\hspace{2mm}}
\newcommand{\Remark}{\note{Remark}}
\newcommand{\forevery}{\quad \forall}

\newcommand{\norm}[1]{\left \|#1\right \|}

\newcommand{\rec}[1]{\frac{1}{#1}}

\newcommand{\dist}{\mathrm{dist}\,}
\newcommand{\sign}{\mathrm{sign}\,}

\newcommand{\prt}{\partial}
\newcommand{\sms}{\setminus}

\newcommand{\ti}{\times}

\newcommand{\tl}{\tilde}

\newcommand{\sbs}{\subset}

\newcommand{\tin}{\to\infty}
\newcommand{\ind}[1]{_{_{#1}}}

\newcommand{\sth}{such that\xspace}

\newcommand{\bvp}{boundary value problem\xspace}

\newcommand{\bdw}{\partial\Gw}

\newcommand{\qtxt}[1]{\quad\textrm{#1}}

\def\ga{\alpha}            \def\gg{\gamma}
       \def\gd{\delta}      \def\ge{\epsilon}

            \def\gl{\lambda}
                 
    \def\gr{\rho}        
\def\gs{\sigma}       
      \def\gw{\omega}
\def\gx{\xi}                \def\gz{\zeta}
     \def\Gd{\Delta}

\def\Gw{\Omega}              

\def\BBG {\mathbb G}       
   \def\BBK {\mathbb K}    
   \def\BBN {\mathbb N}    
   \def\BBR {\mathbb R}

\def\GTM {\mathfrak M}

\def\btr{boundary trace\xspace}
\def\tr{\mathrm{tr}}

\def\LVsup{$L_V$ superharmonic\xspace}
\def\LVsub{$L_V$ subharmonic\xspace}

\def\q{\quad}

\begin{document}

\title[Boundary value problems with signed measure data]{Boundary value problems with signed measure data for semilinear\\ Schr\"odinger equations}

\author{Moshe Marcus}

\date{\today}

\begin{abstract}
		Consider operators  $L_{V}:=\Gd + V$ in a bounded Lipschitz domain $\Gw\sbs \mathbb{R}^N$. Assume that $V\in C^\ga(\Gw)$ satisfies $|V(x)| \leq \bar a\,\dist(x,\bdw)^{-2}$ in $\Gw$ and that $L_V$ has a (minimal) ground state $\Phi_V$ in $\Gw$.
	We derive a representation formula for signed supersolutions (or subsolutions) of $L_Vu=0$ possessing an $L_V$ \btr. We apply this formula to the study of some questions of existence and uniqueness  for an associated semilinear \bvp with signed measure data. \\ [2mm]
		[2mm]
\end{abstract}

\maketitle
\section{Introduction}
	Consider equations of the form
\begin{equation}\label{main_eq}\BAL
-L_Vu+f\circ u=\tau \qtxt{in }\Gw
\EAL\end{equation}
where $L_V=\Gd+V$, $\Gw$ is a bounded Lipschitz domain in $\BBR^N$, $f\in C(\Gw\ti \BBR)$ and $f\circ u(x)=f(x,u(x))$.

We assume that the potential $V\in C^\ga(\Gw)$ for some $\ga>0$ and satisfies the following conditions:
\[ \exists\, \bar a>0\, :\quad |V(x)| \leq \bar a \gd(x)^{-2} \forevery\, x\in \Gw, \tag{A1}  \]
$$\gd(x):=\dist(x,\bdw) $$
and
\[ \inf A< 1 < \sup A, \tag{A2}\]
$$ A:= \{\gg: \int_\Gw |\nabla\phi|^2\,dx\geq \gg\int_\Gw \phi^2 V\, dx  \forevery \phi\in H^1_0(\Gw)\}.$$

These conditions imply the existence of a (minimal) Green function $G_V$ and a Martin kernel $K_V$ for the operator $-L_V$ \cite{An87}. Here, the Martin kernel is a function $K_V$ on $\Gw\ti \bdw$ \sth for every point $y\in \bdw$ the function $\Gw\ni x\mapsto K_V(x,y)$ is a positive $L_V$ harmonic function that vanishes on $\bdw\sms \{y\}$ \cite{An87}. 
Moreover, the ground state $\Phi_V$ is the unique positive eigenfunction of $-L_V$ with eigenvalue $\gl_V>0$. The ground state and the Martin kernel are normalized by setting   $$\Phi_V(x_0) = 1, \q K_V(x_0,y)=1  \forevery y\in \bdw,$$
where $x_0$ is a fixed reference point in $\Gw$.

It is known that, for any $y\in\Gw$ and $\ge>0$, there exist constants $c_i>0$ \sth,
\begin{equation}
c_1\Phi_V(x) \leq G_V(x,y)\leq c_2\Phi_V(x) \qtxt{for }\; x\in \Gw: \, |x-y|>\ge. 
\end{equation} 
This follows, for instance, from \cite[Theorem 3.1 and Lemma 3.6]{Pincho89}. 

Consequently, if $\tau\in \GTM(\Gw,\Phi_V)$ i.e., $\int_\Gw \Phi_V\,d|\tau|<\infty$ then, 
\begin{equation}\label {Gtau}
\BBG_V[|\tau|](x):=\int_\Gw G_V(x,y)d|\tau|(y)< \infty \forevery\, x\in \Gw.
\end{equation}
For any positive measure $\tau$, either \eqref {Gtau} holds or $\BBG_V[\tau]\equiv\infty$ in $\Gw$.

We consider the \bvp,
\begin{equation}\label{bvp-f}
\BAL
-L_Vu+f\circ u&=\tau \qtxt{in }\Gw\\
\tr_V u&=\nu \qtxt{on }\bdw,
\EAL
\end{equation}
where $\tau\in \GTM(\Gw;\Phi_V)$, $\nu\in \GTM(\bdw)$. Here $\tr_V$ denotes the $L_V$ boundary trace defined via harmonic measures for $L_V$ in domains $D\Subset \Gw$, similar to the classical definition of \btr for positive harmonic functions (see Section 2.1 below). 

With respect to the nonlinear term we assume that $f\in C(\Gw\ti\BBR)$ and
\begin{equation}\label{f-cond}\BAL
&f(x,\cdot)  \text{ is non-decreasing} \forevery x\in \Gw,\\ 
&f(\cdot,t)\in L^1(\Gw;\Phi_V) \forevery t\in \BBR.
\EAL\end{equation}

A function $u\in L^1_{loc}(\Gw)$ \sth $f \circ u\in L^1(\Gw;\Phi_V)$ is a solution of \eqref{main_eq} if the equation holds in the distribution sense. The function $u$ is a sub or supersolution of the equation if the appropriate inequlity holds in the distribution sense.
The term $L_V$ harmonic or $L_V$  superharmonic is reserved for the case  $-L_Vu=0$ or $-L_Vu\geq 0$.

 A function $u\in L^1_{loc}(\Gw)$ is a solution of \eqref{bvp-f} if $f \circ u\in L^1(\Gw;\Phi_V)$, the equation holds in the distribution sense and the boundary data is attained 
as an $L_V$ trace. Similarly $u$ is a supersolution of \eqref{bvp-f} if $f \circ u\in L^1(\Gw;\Phi_V)$, $-L_Vu+f \circ u\geq \tau$ in the sense of distributions and $u$ has an $L_V$ trace 
$\nu'\geq \nu$. A subsolution is defined in the same way, with $\geq$ replaced by $\leq$.

A positive \LVsup function $v$ is an $L_V$ \textit{potential} if it does not dominate any positive $L_V$ harmonic function. Under the present conditions $v$ is an $L_V$ potential if and only if \cite {An88}, 
$$ \exists\tau\in \GTM_+(\Gw;\Phi_V): \q v=\BBG_V[\tau].$$
Using this fact and the Martin representation theorem as in \cite{An87}, the Riesz decomposition lemma may be stated as follows:

\begin{lemma}\label{Riesz}
If $V$ satisfies (A1) and (A2), a positive  $L_V$ superharmonic function $w$ can be written in the form
\begin{equation}\label {Riesz}
w=\BBK_V[\nu] + \BBG_V [\tau] \qtxt{where }\; \tau\in \GTM_+(\Gw;\Phi_V), \; \nu\in \GTM_+(\bdw).
\end{equation}
\end{lemma}

In the present paper we start with a brief disscussion of the $L_V$ trace, always assuming (A1) and (A2). We show that (see Section 2.1) 
\[\BAL \tr_V \BBK_V[\nu]&=\nu \forevery \nu\in \GTM(\bdw) \\
\tr_V \BBG_V[\tau] &=0 \forevery \tau\in \GTM(\Gw;\Phi_V).
\EAL \]

Using this result we show that a formula similar to \eqref{Riesz} holds for \textit{signed} $L_V$ sub or superharmonic functions provided that the function has an $L_V$ trace (see Theorem \ref{add4} below).

In Section 3 we study problem \eqref{bvp-f} and prove:\\ [1mm]
 (a)  If $u_1$ is a supersolution and $u_2$ a subsolution of 
 \eqref{bvp-f} then $u_2\leq u_1$. Consequently \textit{\eqref{bvp-f} has at most one solution.}\\ [1mm]
(b) If \eqref{main_eq} has a subsolution $u_2$ and a supersolution $u_1$ \sth the $L_V$ trace $\nu_i:=\tr_Vu_i$ exists and $\nu_2\leq \nu_1$ then problem
\eqref{bvp-f} has a solution for every $\nu\in \GTM(\bdw)$ \sth $\nu_2\leq \nu\leq \nu_1$.

Results similar to (b) involving \textit{weak} sub and supersolutions of the \bvp \eqref{bvp-f} with $V=0$ have been obtained in \cite{MV-book} and \cite{Mon-Po}. The present interpretation of boundary data (as $L_V$ trace) is, of course, different from its interpretation in the weak formulation.

For other related results see  \cite{CS87} and  \cite{DH74}.

\section{Boundary trace and a representation formula}
\subsection{$L_V$ trace}

Let $D\Subset \Gw$ be a Lipschitz domain and denote by $\gw_V^{x_0,D}$ the harmonic measure of $L_V$
in $D$ relative to a point $x_0\in D$.
If $P_V^D$ is the Poisson kernel of $L_V$ in $D$ then
\begin{equation}\label{hr-m1}
d \gw_V^{x_0,D} = P_V^D(x_0,\cdot)dS \qtxt{on } \; \prt D.
\end{equation}

Let $\{D_n\}$ be a uniformly Lipschitz exhaustion of $\Gw$. It is well known that if  $u$ is a positive $\Gd$-harmonic function then
\begin{equation}\label{gw0}
u\lfloor_{_{\prt D_n}} d\gw_0^{x_0, D_n} \rightharpoonup \nu
\end{equation}
where $\nu\in \GTM(\bdw)$ is the boundary trace of $u$ and $\rightharpoonup$ indicates weak convergence in measure. Similarly we define,

\begin{definition}\label{ntr}  	A non-negative Borel function $u$ defined in $\Gw$ has an $L_V$ \emph{boundary trace} $\nu\in \GTM(\bdw)$ if
	\begin{equation}\label{ntr1}
	\lim_{n\to\infty}\int_{\prt D_n}hud\gw_V^{x_0,D_n}= \int_{\bdw} hd\nu \forevery h\in C(\bar\Gw),
	\end{equation}
	for every uniformly Lipschitz exhaustion $\{D_n\}$ of $\Gw$.
	The $L^V$ trace will be denoted by $\tr_V(u)$.
	(Here  we assume that the reference point $x_0$ is in $D_1$.) The $L_V$ \btr of $u$ is denoted by $\tr_Vu$.
	
	A signed Borel function $u$ has an $L_V$ \btr $\nu$ if
	$$\sup \int_{\prt D_n}|u|d\gw_V^{x_0,D_n} <\infty$$
and \eqref{ntr1} holds.	
\end{definition}

\begin{lemma}\label{weak_con} Let $\{D_n\}$ be a uniformly Lipschitz exhaustion of $\Gw$. Assume that $K_V$ is normalized at $x_0\in \Gw$. Then, for every positive $L_V$ harmonic function $u=\BBK_V[\nu]$,
	\begin{equation}\label{hr-meas}
	\lim_{n\to \infty}\int_{\prt D_n}hu\, d\gw_V^{x_0,D_n} = \int_{\bdw}h\,d\nu
	\forevery h\in C(\bar \Gw).
	\end{equation}
\end{lemma}

\Remark The proof is similar to that of \cite[Lemma 2.2] {MV-Lip}. For the convenience of the reader, we provide it below.

\proof First observe that, by \eqref{hr-m1},
\begin{equation}\label{hr-m2}
u(x_0)= \int_{\prt D_n}u\, d\gw_V^{x_0,D_n}
\end{equation}
and, as the Martin kernel is normalized by $K_V(x_0,y)=1$ for every $y\in \bdw$,
$$u(x_0)=\int_{\bdw}K_V(x_0,y)\,d\nu(y)= \nu(\bdw). $$
Thus \eqref{hr-meas} holds for $h\equiv 1$. By \eqref{hr-m2} the following sequence of measures is bounded:
$$\gs_n=\begin{cases}u\, d\gw_V^{x_0,D_n} & \text{on }\;\prt D_n\\
0 & \text{on }\;\bar\Gw\sms\prt D_{n}
\end{cases}  \q  n\in \BBN. $$
Let $\{\gs_{n_k}\}$ be a weakly convergent subsequence with limit $\nu'$. Then $$\nu(\bdw)=\nu'(\bdw)= u(x_0).$$

Let $F\sbs \bdw$ be a compact set and define
$$\nu^F= \nu\chr{F}, \q  u^F=\BBK_V[\nu^F].$$
Let $\gs_n^F$ be defined in the same way as $\gs_n$ with $u$ replaced by $u^F$.  Proceeding as before  we obtain a weakly convergent subsequence of $\{\gs_n^F\}$
with limit $\hat \nu_F$ supported in $F$. Furthermore,
$$\hat \nu_F(F)=\hat \nu_F(\bdw)=u^F(x_0)=\nu^F(\bdw)=\nu(F).$$
Since $u^F\leq u$ it follows that $\hat \nu_F\leq \nu'$. Consequently
$\nu(F) \leq \nu'(F)$. As this inequality holds for every compact subset of $\bdw$ it follows that $\nu\leq \nu'$. As the measures are positive and $\nu(\bdw)=\nu'(\bdw)$ it follows that $\nu=\nu'$ and therefore the whole sequence $\{\gs_n  \}$ converges weakly to $\nu$.

\qed

\begin{lemma}\label{trKV}
	Assume (A1) and (A2). Then
	\begin{equation}\label{trKV1}\BAL
	(a)\q	&\tr_V(\BBK_V[\nu])=\nu &&\forevery \nu\in \GTM(\bdw)\\
	(b)\q 	&\tr_V(\BBG_V[\tau])=0  &&\forevery \tau\in \GTM(\Gw;\Phi_V).
	\EAL\end{equation}
\end{lemma}

\proof
(a) This is proved in Lemma \ref{weak_con} for $\nu\geq 0$ and, by linearity, it holds for every $\nu\in \GTM(\bdw)$.   \\ [2mm]
(b) Without loss of generality we may assume $\tau\geq 0$. Then $u:=\BBG_V[\tau]$ is an $L_V$ potential , i.e., a positive $L_V$ superharmonic function that does not dominate any positive $L_V$ harmonic function (see \cite{An88}). 

 Let $\{D_n\}$ be a smooth exhaustion of $\Gw$ and let
\begin{equation}\label{u_in_Dn}
u_n(x)= \int_{\prt D_n}P_V^{D_n}(x,\gx)u(\gx)dS_\gx \forevery x\in D_n, 
\end{equation}
where $u\lfloor_{\prt D_n}$ is the trace (in the Sobolev sense) of $u$ on $\bdw_n$.
Then $u_n$ is $L_V$ harmonic, $u_n\leq u$ and $\{u_n\}$ is decreasing. Therefore the limit $\underline{u}:=\lim u_n$ is non-negative, $L_V$ harmonic and $\underline{u}\leq u$. Since $u$ is an $L_V$ potential it follows that $\underline{u}=0$. Thus $\tr_V u=0$.
\qed

\begin{corollary}\label{trVu} 
  If $u=K_V[\nu] + G_V[\tau]$ for $\nu$  and $\tau$ as in \eqref{trKV1} and $\{D_n\}$ is a smooth exhaustion of $\Gw$ then,
\begin{equation}\label{hr-gen2}
\int_{\prt D_n}u(\gx)\, d\gw_V^{x,D_n}(\gx) \to K_V[\nu] \forevery x\in \Gw.
\end{equation}
\end{corollary} 

\proof
If $u=K_V[\nu]$ then,
\begin{equation}\label{hr-gen1}
u(x)= \int_{\prt D_n}u(\gx)\, d\gw_V^{x,D_n}(\gx)= \int_{\prt D_n}u(\gx) P_V^{D_n}(x,\gx) dS_\gx  
\end{equation}
for every $x\in D_n$. Therefore \eqref{hr-gen2} holds.

 If $u=G_V[\tau]$, $\tau\in \GTM_+(\Gw;\Phi_V)$ then, by  Lemma \ref{trKV},
\begin{equation}\label{Gtau1}
 \int_{\prt D_n}G_V[\tau]\, d\gw_V^{x,D_n}(\gx)\to 0 \forevery x\in \Gw. 
\end{equation}
(For given $x\in \Gw$ the limit is taken over sufficiently large $n$ \sth $x\in D_n$.)
Therefore, by linearity, \eqref{Gtau1} holds for every $\tau\in \GTM(\Gw;\Phi_V)$. 

 \qed

\subsection{A representation formula.}

\begin{lemma}\label{chk1}
	Let $w\in L^1_{loc}(\Gw)$ be an \LVsub function, i.e. $-L_Vw\leq 0$, in the sense of distributions. 

Suppose that there exists a smooth exhaustion $\{\Gw_n\}$ \sth 
\begin{equation}\label{sup_finite}
\sup	\int_{\bdw_n}P_{V,n}(x_0,y)w(y) dS <\infty.
\end{equation}

Then 
\begin{equation}\label{eq:add0}
\gl:=L_Vw \in \GTM_+(\Gw;\Phi_V).
\end{equation}
Thus $v:= w + \BBG_V[\gl]$ is $L_V$ harmonic and
\begin{equation}\label{v=lim}
v(x)=\lim_{n\tin} \int_{\bdw_n}P_{V,n}(x,y)w(y) dS, \forevery x\in \Gw.
\end{equation}
If $w\geq 0$,
\begin{equation}\label{add0'}
\exists \nu\in \GTM_+(\bdw) \Longrightarrow w+\BBG_V[\gl]=\BBK_V[\nu].
\end{equation}

\end{lemma}

\proof There exists a positive Radon measure $\gl$ \sth
$$-L_Vw=-\gl.$$
Let $D\Subset\Gw$ be a Lipshitz domain. Since $\gl\chr{D}$ is a finite measure and $Vw\chr{D}\in L^1(D)$, $w+ G_0^D((\gl-Vw)\chr{D})$ is a harmonic funtion in $D$. Hence $w\in W^{1,p}_{loc}(D)$ for $1\leq p<N/(N-1)$. 
Consequently, as a Sobolev trace, $w\lfloor_{\prt D}\in L^p(\prt D)$.

Let $\{\Gw_n\}$ be the smooth exhaustion of $\Gw$ \sth \eqref{sup_finite} holds. Put
$ \gl_n=\gl\chr{\Gw_n}$ and let $h_n$ be the Sobolev trace of $w$ on $\bdw_n$.

Then $w + \BBG_{V,n}[\gl_n]$ is $L_V$ harmonic in $\Gw_n$ and 
\begin{equation}\label{gln}
w + \BBG_{V,n}[\gl_n]= \int_{\bdw_n}P_{V,n}(\cdot,y)h_n(y) dS \qtxt{in } \Gw_n.
\end{equation}
Hence, by \eqref{sup_finite}, the sequence $\{\BBG_{V,n}[\gl_n](x_0)\}$ - which is non-decreasing - is bounded above and, by the monotone convergence theorem,
$$ \BBG_{V,n}[\gl_n](x_0) \to  \BBG_V[\gl](x_0)<\infty.$$
This implies that $\gl\in \GTM_+(\Gw;\Phi_V)$ and that 
$$ w + \BBG_{V,n}[\gl_n] \to w + \BBG_V[\gl] =: v$$
everywhere in $\Gw$. It also implies that the right hand side of \eqref{gln} converges and \eqref{v=lim} holds.

The last assertion follows from the Martin representation theorem (Ancona \cite{An87}).

\qed

Here is a simple consequence of the lemma:
\begin{corollary}\label{c:check1}
	Let $w\in L^1_{loc}(\Gw)$ be a non-negative \LVsub function.
Suppose that there exists a smooth exhaustion $\{\Gw_n\}$ \sth 
\begin{equation}\label{sup_finite_0}
\lim_{n\to\infty}	\int_{\bdw_n}P_{V,n}(x_0,y)w(y) dS = 0.
\end{equation}
Then $w\equiv 0$.
\end{corollary}
\proof The function $w$ satisfies the assumptions of Lemma \ref{chk1}. Let $v$ be as in the statement of the lemma. Then $v$ is $L_V$ harmonic and $v\geq 0$. By \eqref {v=lim} and \eqref{sup_finite_0},  $v(x_0)= 0$. By Harnack's inequality, $v\equiv 0$ and (as $\gl\geq 0 $) $w\equiv 0$.
\qed 

\begin{lemma}\label{chk2}
	Suppose that $u\in L^1_{loc}(\Gw)$ is $L_V$ harmonic and that 
	\begin{equation}\label{sup-abs}
	\sup	\int_{\bdw_n}P_{V,n}(x_0,y)|u(y)| dS <\infty, 
	\end{equation}
	for some smooth exhaustion $\{\Gw_n\}$ of $\Gw$.
	
	Then  $\tr_Vu$ exists and  $u= \BBK_V[\nu]$ for some $\nu\in \GTM(\bdw)$.
\end{lemma}

\proof By Kato's inequality $-L_Vu_+\leq 0$. Assumption \eqref{sup-abs} implies that $u_+$ satisfies \eqref{sup_finite}. Therefore, by Lemma \ref{chk1}, 
$$\gl':= L_Vu_+ \in \GTM_+(\Gw;\Phi_V).$$
Thus $u_+ + \BBG_V[\gl']$ is a positive $L_V$ harmonic function. By the reprentation theorem, there exists $\nu'\in \GTM_+(\bdw)$ \sth,
\begin{equation}
u_+ + \BBG_V[\gl'] = K_V[\nu'].
\end{equation}
Applying the same argument to $-u$ we conclude that there exist 
$\gl''\in \GTM_+(\Gw;\Phi_V)$ and $\nu''\in \GTM_+(\bdw)$ \sth
\begin{equation}
u_- + \BBG_V[\gl''] = K_V[\nu''].
\end{equation}

Hence,
$$u= \BBG_V[\gl'' - \gl'] + K_V[\nu'-\nu''].$$
Since $u$ is $L_V$ harmonic, $\gl'' -\gl'=0$. Thus $u=K_V[\nu'-\nu'']$ and $\tr_V u=\nu'-\nu''$.

\qed

The following is an immediate consequence.
\begin{corollary}\label{tr-exist'}
If $u\in L^1_{loc}(\Gw)$, $-L_Vu=:\tau\in \GTM(\Gw;\Phi_V)$ and \eqref{sup-abs}
holds then $u$ has an $L_V$ \btr, say $\nu$, and $u=\BBG_V[\tau]+\BBK_V[\nu]$.
\end{corollary}

Using the above lemmas we obtain,

\begin{theorem}\label{add4} Let $w\in L^1_{loc}(\Gw)$.
Suppose that $-L_Vw \leq \tau$ (resp. $-L_Vw \geq \tau$ ) for some  $\tau \in \GTM(\Gw;\Phi_V)$ and
 \begin{equation}\label{sup|w|}
\sup \int_{\prt \Gw_n}|w|d\gw_V^{x_0,\Gw_n} <\infty.
 \end{equation}
Then there exists $\gl\in \GTM_+(\Gw;\Phi_V)$ \sth 
\begin{equation}\label{LVw}
-L_Vw =\tau-\gl,  \q (\mathrm{resp. } -L_Vw=\tau+\gl).
\end{equation}

Moreover $w_+$ and $w_-$ have an $L_V$ trace; thus $\tr_Vw$ exists. If $\gs:= \tr_Vw$ then 
\begin{equation}\label{eq:add'}
w = \BBG_V[\tau-\gl] + \BBK_V[\gs] \q (\mathrm{resp. }\; w = \BBG_V[\tau+\gl] + \BBK_V[\gs]).
\end{equation}
\end{theorem}

\proof  It is sufficient to prove the result in the case $-L_Vw\leq \tau$. If $-L_Vw \geq \tau$ the problem is reduced to the previous case by considering $-L_V(-w)\leq -\tau$.

If $-L_Vw \leq \tau$ then $w^*=w-\BBG_V[\tau]$ is $L_V$ subharmonic. Moreover, as $\tr_V \BBG_V[\tau]=0$ and $w$ satisfies \eqref{sup|w|},  that $w^*$ satisfies \eqref{sup|w|}. Therefore, by Lemma \ref{chk1}, there exists $\gl\in \GTM_+(\Gw;\Phi_V)$ \sth
\begin{equation}\label{v-har}
v:= w^* +\BBG_V[\gl] \qtxt{is $L_V$ harmonic.}
\end{equation} 
This implies \eqref{LVw}.
Furthermore, by Lemma \ref{chk2}, $\gs:=\tr_V v$ exists and $v=\BBK_V[\gs]$.
This implies \eqref{eq:add'}. 

By an extension of Kato's inequality \cite{BP'} (see also cite \cite[Corollary 1.5.7]{MV-book}), $-L_Vw_+ \leq \tau_+$. Obviously $w_+$ satisfies \eqref{sup|w|}. Therefore, by the first part of the proof,
$\tr_V w_+$ exists. 

Further, by \eqref{LVw}, $-L_V(-w)=\gl-\tau\in \GTM(\Gw;\Phi_V)$. Therefore $w_-=(-w)_+$ has an $L_V$ trace.

\qed

To simplify statements we introduce the following definition:

A function $u$ is \textit{$L_V$ perfect} if
\begin{equation}\label{e:perfect}
u=\BBG_V[\tau] + \BBK_V[\nu] \qtxt{where }\;\tau\in \GTM(\Gw;\Phi_V),\;\nu\in \GTM(\bdw).
\end{equation} 
 
\remark\label{r:perfect} Obviously an $L_V$ perfect function has an $L_V$ trace. By Theorem \ref{add4}, if $w$ is \LVsub or \LVsup 
 $$\tr_V w \;\text{exists}\q \Longleftrightarrow \qtxt{$w$ is $L_V$ perfect.}$$

The following is a simple consequence of the Riesz decomposition lemma.
\begin{lemma}\label{tr_exist}
	Let $u\in L^1_{loc}(\Gw)$ be a positive function \sth $-L_Vu=:\tau\in \GTM(\Gw;\Phi_V)$. Then $u$ is $L_V$ perfect.
\end{lemma}

\proof
Put  $U=u+ G_V[|\tau|-\tau]$. Then $U$ is positive, \LVsup and $-L_VU=|\tau|$. By the Riesz decomposition lemma, there exists 
$\nu\in \GTM_+(\bdw)$ \sth 
$U=\BBG_V[|\tau|] +\BBK_V[\nu]$. Thus 
$$u=U-G_V[|\tau|-\tau]=\BBG_V[\tau] + \BBK_V[\nu].$$
By Lemma \ref{trKV}, $\tr_V u=\nu$.
\qed

\vskip 4mm

\section{The nonlinear problem}
We consider equation \eqref{main_eq} and the \bvp \eqref{bvp-f}
where 
$$ \tau\in \GTM(\Gw;\Phi_V), \q \nu \in \GTM(\partial \Omega).$$
As before we assume that $V$ satisfies conditions (A1), (A2) and that $f\in C(\Gw;\BBR)$ satisfies condition \eqref{f-cond}.

\subsection{$L_V$ perfect functions.}

\begin{lemma}\label{l:subsup_rep}
	Assume that $u\in L^1_{loc}(\Gw)$ and $f\circ u\in L^1(\Gw;\Phi_V)$.
	
If $u$ is a subsolution (resp. supersolution) of equation \eqref{main_eq} and $\tr_Vu$ exists then $u$ is $L_V$ perfect.
More precisely, there exists a measure $\gl\in \GTM_+(\Gw;\Phi_V)$ \sth, if $\gs:=\tr_V u$ then,
\begin{equation}\label{e:subsup_rep}\BAL
u &=  G_V[\tau-f\circ u-\gl]+K_V[\gs] \\ (u &=  G_V[\tau-f\circ u+\gl]+K_V[\gs]). 
\EAL\end{equation}
\end{lemma}

\proof If $u$ is a subsolution of \eqref{main_eq} then $-L_Vu \leq \tau-f\circ u\in \GTM(\Gw;\Phi_V)$.
 Therefore, by Theorem \ref {add4}, $u$ is $L_V$ perfect and satisfies \eqref{e:subsup_rep}. A similar argument applies to supersolutions.
 
\qed

\begin{lemma}\label{LVsup1} 
 Let $w\in L^1_{loc}(\Gw)$ be a positive function \sth $f \circ w\in L^1(\Gw;\Phi_V)$ and 
 $$-L_Vw+f \circ w\geq0.$$
Then $w$ has an $L_V$ boundary trace, say $\nu$, and
$$-L_V w + f \circ w =:\gl\in \GTM_+(\Gw;\Phi_V).$$
Thus
	\begin{equation}\label{sup-rep} 
w= G_V[\gl-f \circ w] + K_V[\nu].
\end{equation}
\end{lemma}

\proof Apply the Riesz decomposition lemma to $w+G_V[f \circ w]$ which is a positive \LVsup function. 
 
  \qed
 
\begin{lemma}\label{tru<0}
Suppose that $u$ is an $L_V$ perfect function.

(i) If $\tr_Vu\leq 0$ then $\tr_Vu_+=0$,

(ii) $|u|$ is dominated by an $L_V$ superharmonic function.
\end{lemma}

\proof By assumption, $u$ satisfies \eqref{e:perfect}:
$u=G_V[\tau] + K_V[\nu]$. Hence if $\nu\leq 0$ then $u_+\leq G_V[\tau_+]$  and, by Lemma \ref {trKV}, $\tr_V u_+=0$. 

Moreover,
 $$|u|\leq G_V[|\tau|] +K_V[|\nu|]$$
 and the right hand side is \LVsup.
\qed

\begin{lemma}\label{u1<u2}
(i) Let $u_1$ (resp. $u_2$) be a supersolution (resp. subsolution)
of \eqref{main_eq} and assume that $f(u_i)\in L^1(\Gw;\Phi_V)$ and $\tr_Vu_i$ exists $(i=1,2)$. If $\tr_Vu_2\leq \tr_Vu_1$ then $u_2\leq u_1$.\\ [1mm]
(ii) Problem  \eqref{bvp-f} has at most one solution.
\end{lemma}

\proof 
(i) The assumptions imply that 
$$-L_V (u_2 - u_1) + (f\circ u_2- f\circ u_1)\leq 0.$$
Therefore, by Kato's lemma (see e.g.  \cite[Proposition 1.5.4]{MV-book}),
$$ -L_V(u_2-u_1)_+ + \sign_+(u_2-u_1) (f\circ u_2- f\circ u_1) \leq 0.$$
Since $f$ is non-decreasing, the second term is non-negative. Hence $(u_2-u_1)_+$ is $L_V$ subharmonic.

By Lemma \ref{l:subsup_rep}, $u_1$ and $u_2$ are $L_V$ perfect. Let $u:=u_2-u_1$. As $\tr_Vu\leq 0$, Lemma \ref{tru<0}(i) implies that $\tr_Vu_+=0$. Hence, by Corollary \ref{c:check1}, $u_+=0$. Thus $u_2\leq u_1$.

(ii) is an immediate consequence of (i).

\qed
\subsection{The sub-supersolution method for problem \eqref{bvp-f}.}

\begin{lemma}\label{GV-W}
The operator $\BBG_V$ maps $\GTM(\Gw;\Phi_V)$ into  $W^{1,p}_{loc}(\Gw)$, $1\leq p<N/(N-1)$. For every domain $D\Subset \Gw$, i.e. $\bar D\subset \Gw$,
\begin{equation}\label{GV-bd}
\norm{\BBG_V[\tau]}^p_{W^{1,p}(D)}\leq c_1(D,p) \int_\Gw\Phi^p_V d|\tau|.
\end{equation}
In particular $\BBG_V:\GTM(\Gw;\Phi_V) \to W^{1,1}_{loc}(\Gw)$ is a bounded mapping.
\end{lemma}

\proof Let $D$ be a Lipschitz domain \sth $\bar D \subset \Gw$. Given $\tau\in \GTM(\Gw;\Phi_V)$ put
$\tau_ D=\tau\chr{D}$. Then,
\begin{equation}\label{CGV1}
\norm{\BBG_V[\tau_D]}^p_{W^{1,p}(D)}\leq c_0(D,p) \norm{\tau_D}_{\GTM(\Gw)} 
\leq c_1(D,p) \int_\Gw\Phi_V(y)d|\tau_D|.
\end{equation} 
 The first inequality in \eqref{CGV1} is a consequence
of the following: by \cite{Hrv}, 
$$G_V(x,y) \leq C_1(D) |x-y|^{2-N} \forevery (x,y)\in D, \;x\neq y.$$
Hence, by interior elliptic estimates in $D\sms \{y\}$,
$$|\nabla_x G_V(x,y)|\leq C'_1(D) |x-y|^{1-N}
\forevery x,y\in D,\; x\neq y.$$

The second inequality in \eqref{CGV1} is a consequence of  Harnack's inequality and the normaliztion $\Phi_V(x_0)=1$. These imply that $\Phi_V$ is bounded away from zero in every compact subset of $\Gw$.

Let $\rho_D = \dist(D,\bdw)$ and put 
$$E=\{y: \gd(y)\leq \rho'_D\}, \q  \rho'_D:= \rho_D/ 16(1+\kappa)^2$$
where $\kappa\geq 1$ denotes a Lipschitz constant for $\bdw$.
By \cite[Theorem 1.4]{MM-Green},
$$G_V(x,y) \leq C \frac {\Phi_V(x)\Phi_V(y)}{\Phi_V^2(x_y)}|x-y|^{2-N}$$
for every $x\in D$, $y\in E$. Here
$x_y$ is a point \sth $\gd(x_y)\geq |x-y|$. As $|x-y|>\rec{2}\gr_D$,
\begin{equation}
G_V(x,y)\leq c_2(D) \Phi_V(y) \forevery (x,y)\in D\ti E.
\end{equation}
Hence, by interior elliptic estimates,
$$\sup_{x\in D}|\nabla_x G_V(x,y)|\leq c_2'(D)\Phi_V(y) \forevery y\in E.$$

Consequently, if $\tau_E := \tau\chr{E}$, 
\begin{equation}
\norm{\BBG_V[\tau_E]}_{W^{1,p}(D)}^p \leq c_2(D,p) \int_\Gw\Phi^p_V(y)d|\tau_E|.
\end{equation}
This inequality together with \eqref{CGV1} applied to $\tau\chr{\Gw\sms E}$, imply the stated result.

\qed

\begin{theorem}\label{exists1}
Let $u_1$ be a supersolution and $u_2$ a subsolution of equation \eqref{main_eq}. Suppose that $f(u_i)\in L^1(\Gw;\Phi_V)$ and that $u_i$ satisfies \eqref{sup|w|} $(i=1,2)$. If 
\begin{equation}\label{temp}
u_2+ \BBG_V[f\circ u_2]\leq u_1 + \BBG_V[f\circ u_1]
\end{equation}
then there exists a solution $u$ of 
\eqref{bvp-f} \sth $u_2\leq u\leq u_1$.

Further, $\nu_i:=\tr_V u_i$, $i=1,2$, exist and $\nu_2\leq \nu_1$. Finally, for every $\nu\in \GTM(\bdw)$ \sth $\nu_2\leq\nu \leq \nu_1$, \bvp \eqref{bvp-f} has a unique solution.
\end{theorem}

\Remark For similar results regarding \textbf{weak} sub - supersolutions of \eqref{main_eq} with weak trace and $V=0$ see \cite{MV-book}, \cite{Mon-Po} and the references therein.

\proof By Theorem \ref{add4} applied to $w=u_i + G_V[f\circ u_i]$, there exist measures $\gl_i\in \GTM_+(\Gw;\Phi_V)$ and $\nu_i\in \GTM(\bdw)$ $i=1,2$, \sth

\begin{equation}\label{u1u2}
\BAL
u_1 + G_V[f\circ u_1] &= G_V[\tau +\gl_1] +K_V[\nu_1],\\ 
u_2 + G_V[f\circ u_2] &= G_V[\tau -\gl_2] +K_V[\nu_2].
\EAL
\end{equation}
In particular $\nu_i=\tr_Vu_i$.
Hence, using \eqref{temp},
$$ K_V[\nu_2] - G_V[\gl_2] \leq K_V[\nu_1] + G_V[\gl_1].$$
It follows that,
$$\nu_2 = \tr_V ( K_V[\nu_2] - G_V[\gl_2])\leq \nu_1= \tr_V (K_V[\nu_1] + G_V[\gl_1]).$$
Therefore, by Lemma \ref{u1<u2}, $u_2\leq u_1$.

Problem \eqref{bvp-f} is equivalent to the following: find $u\in L^1_{loc}(\Gw)$ \sth $f\circ u\in L^1(\Gw;\Phi_V)$ and
\begin{equation}\label{bvp_f_alt}
u=\BBG_V[\tau-f\circ u] +\BBK_V[\nu].
\end{equation} 
Let $\tl f$ be the function given by
\begin{equation}\label{tlf1}
\tl f(x,t)=\begin{cases} f(x,t) &\text{if }u_2(x)\leq t\leq u_1(x),\\ f(x,u_1(x)) &\text{if }u_1(x)\leq t,\\
f(x,u_2(x)) &\text{if }t\leq u_2(x).\end{cases}
\end{equation}
Note that,
\begin{equation}\label{tlf2}
|\tl f|\leq |f\circ u_1| + |f\circ u_2|.
\end{equation}
Let $z=T\gz$ , $\gz\in L^1_{loc}(\Gw)$, be the mapping given by 
\begin{equation}\label{temp3}
z= \BBG_{V}[\tau- \tl f\circ \gz] + \BBK_V[\nu]. 
\end{equation}

Assume that $\nu_2\leq \nu\leq \nu_1$. Then,
in view of \eqref{tlf2} and Lemma \ref{GV-W}, for every domain $D\Subset\Gw$,
$$\norm{z}_{W^{1,1} (D)}\leq \int_\Gw \Phi_V d|\tau| + \sum_{i=1}^2 \big (\int_\Gw |f\circ u_i|\, \Phi_V dx  +\norm{\BBK_V[\nu_i]}_{W^{1,1}(D)]}\big ).$$
Thus $T(L^1_{loc}(\Gw))$ is bounded in $W^{1,1}(D)$. Therefore, by the Rellich-Kondrachev imbedding theorem, $T(L^1_{loc}(\Gw))$ is precompact in $L^1(D)$. Thus $T$ maps $L^1_{loc}(\Gw)$ into a compact subset of this space. By the Schauder -- Tychonoff fixed point theorem (see e.g. \cite{GD-fpt}), T has a fixed point $u\in L^1_{loc}(\Gw)$, i.e.
\begin{equation}\label{bvp-tlf}
u=\BBG_V[\tau-\tl f\circ u] +\BBK_V[\nu].
\end{equation}
 Since $f\circ u_i= \tl f\circ u_i$, $u_1$ is a supersolution and $u_2$ is a subsolution of problem \eqref{bvp-tlf}. Therefore, by Lemma \ref{u1<u2},  $u_2\leq u \leq u_1$. Hence $f\circ u=\tl f \circ u$ so that $u$ is a solution of \eqref{bvp-f}.

\qed

\begin{corollary}\label {c:exists}
Suppose that $u_1$, $u_2$ are respectively a supersolution and a subsolution of \eqref{main_eq} and that the \btr $\nu_i:= \tr_V u_i$  exists $i=1,2$. If $\nu_2 \leq  \nu_1$ then problem \eqref{bvp-f} has a (unique) solution for every measure $\nu$ \sth $\nu_2\leq\nu\leq \nu_1$. 
\end{corollary}

\proof By Lemma \ref{l:subsup_rep} the assumptions imply that  \eqref{temp} holds. Thus the corollary is an immediate consequence of the theorem.

\qed

\vskip 20mm

\small{Address: Department of Mathematics, Technion, Haifa, Israel
	
	E-mail: marcusm@math.technion.ac.il}

\newpage


\begin{thebibliography}{99}
\bibitem{An87} A. Ancona, \textit{Negatively curved manifolds, elliptic operators and the Martin boundary,} Annals of Mathematics, Second Series,  {\bf 125} (1987), 495-536.
	
\bibitem{An88} A. Ancona, \emph{Th\'eorie du potentiel sur les graphes et les vari´et´es},  in Springer Lecture Notes in Math., vol.1427, (ed. P.L. Hennequin) (1988)  pp.1-–112.


%
%
\bibitem{BP'} H. Brezis and A. C. Ponce, {\em Kato's inequality when $\Delta u$ is a measure}, C. R. Math. Acad. Sci. Paris,
Ser. I 338 (2004) 599-–604.

\bibitem{CS87} P. Cl\'ement and G. Sweers, {\em Getting a solution between sub- and supersolutions without monotone
iteration}, Rend. Istit. Mat. Univ. Trieste 19 (1987), 189--194. 

\bibitem{DH74} J. Deuel and P. Hess, {\em A criterion for the existence of solutions of non-linear elliptic boundary
value problems}, Proc. Roy. Soc. Edinburgh Sect. A 74A (1974/75), 49-–54.

%
%
%
%
%
%
%
%
%
%
%
%
\bibitem{Hrv} Herv\'e, R.-M., \emph{Recherches axiomatiques sur la thorie des fonctions surharmoniques
et du potentiel,} (French) Ann. Inst. Fourier (Grenoble) {\bf 12} (1962), 415--571.
\bibitem{GD-fpt} A. Granas and J. Dugundji \textit{Fixed point theory},
Springer Science+Business Media, New York (2003).

%
\bibitem{MM-Green} M. Marcus, {\em Estimates of Green and Martin kernels for Schr\'odinger operators with singular
potential in Lipschitz domains,} Ann. I. H. Poincar´e AN 36 (2019) 1183–-1200.

%
%
%
%
%
	\bibitem{MV-Lip} M. Marcus and L. Veron, \emph{Boundary trace of positive solutions	of semilinear elliptic equations
	in Lipschitz domains: the subcritical case}, Ann. Scuola Norm. Sup. Pisa Cl. Sci. (5)
Vol. X (2011), 1--73.

\bibitem{MV-book}{M. Marcus and L. V\'{e}ron}, {\em Nonlinear second order elliptic equations involving measures}. De Gruyter Series in Nonlinear Analysis and Applications, 21. De Gruyter, Berlin, 2014.

\bibitem{Mon-Po} M. Montenegro and A. Ponce, \textit{The sub-supersolution method for weak solutions}, Proc. AMS Vol. 136, No. 7 (2008) 2429 -- 2438.





	
%
%

%
%
%
%
%
%
%
%
%
%
%
%
%
%
%
%
%
%
%

%

\bibitem{Pincho89} Y. Pinchover, \emph{Criticality and ground states for second-order elliptic equations,} J. Differ. Eqns. 80 (1989) 237--250.
	
%
%
	
	
\end{thebibliography}
\end{document}